\parindent5mm
\magnification=1200
\baselineskip=14pt
\hsize=160truemm
\vsize=220truemm
\hoffset=4truemm
\voffset=5truemm
\hfuzz=10pt
\headline={\hfil}
\footline={\hfil\tenrm\folio\hfil}


\font\obf=cmbx9 scaled\magstep1

\centerline{\sevenrm Journal of Pure and Applied Algebra 137 (2) (1999)
1-14}
\null\vskip3truecm

\centerline {\obf The dimension of tensor products of k-algebras arising
from pullbacks }
\bigskip

\centerline {\bf S. Bouchiba$^a$ , F. Girolami$^b$ , S. Kabbaj$^c$}
\medskip

\centerline {\sevenrm $^a$Department of Mathematics,  Faculty of Sciences,
University of Mekn\`es, Mekn\`es, Morocco}

\centerline {\sevenrm $^b$Dipartimento di Matematica, Universit\`a degli
Studi
Roma Tre, 00146 Roma, Italy}

\centerline {\sevenrm $^c$Department of Mathematics,  Faculty of Sciences I,
University of F\`es, F\`es, Morocco}
\bigskip

\noindent
----------------------------------------------------------------------------
--------------------------------------
\noindent {\obf Abstract} \medskip
{\sevenrm The purpose of this paper is to compute the Krull dimension of
tensor products of $k-$algebras arising from
pullbacks. We also state a formula for the valuative dimension.}

\noindent
----------------------------------------------------------------------------
--------------------------------------

\vskip1truecm
\noindent {\obf 0. Introduction}
\bigskip

All rings and algebras considered in this paper are commutative with
identity
elements and, unless otherwise specified, are to be assumed to be
non-trivial. All ring-
homomorphisms are unital. Let $k$ be a field. We denote the class of
commutative $k-$algebras
with finite transcendence degree over $k$ by $C$. Also, we shall use
t.d.($A$) to denote the
transcendence degree of a $k-$algebra $A$ over $k$, $A[n]$ to denote the
polynomial ring
$A[X_1, . . ., X_n]$, and $p[n]$ to denote the prime ideal $p[X_1, . . .,
X_n]$ of
$A[n]$, where $p$ is a prime ideal of $A$. Recall that an integral domain
$R$ of finite (Krull)
dimension $n$ is a Jaffard domain if its valuative  dimension, dim$_v(R)$,
is also $n$. Pr\"ufer
domains and noetherian domains are Jaffard domains. We assume familiarity
with this concept, as in
[1], [6] and [10]. Suitable background on  pullbacks is [4], [11], [12] and
[16]. Any unreferenced
material is standard, as in [12] and [17].
\bigskip

In [20] Sharp proved that if $K$ and $L$ are two extension fields of $k$,
then
dim$(K \otimes_k L)$ = min(t.d.($K$), t.d.($L$)). This result provided a
natural starting
point to investigate dimensions of tensor products of somewhat general
$k-$algebras. This was
concretized by Wadsworth in [21], where the result of Sharp was extended to
AF-domains, that is,
integral domains $A$ such that ht($p$) + t.d.($A/p$) = t.d.($A$), for all
prime ideals $p$ of $A$.
He showed that if $A_1$ and $A_2$ are AF-domains, then dim$(A_1 \otimes_k
A_2)$ =
min(dim($A_1$) + t.d.($A_2$), dim($A_2$) + t.d.($A_1$)). He also stated a
formula for
dim$(A \otimes_k R)$ which holds for an AF-domain $A$, with no restriction
on $R$. We recall, at
this point, that an AF-domain is a (locally) Jaffard domain [13].
\bigskip

In [5] we were concerned with AF-rings. A $k-$algebra $A$ is said to be
an AF-ring provided ht($p$) + t.d.($A/p$) = t.d.($A_p$), for all prime
ideals $p$ of $A$ (for nondomains, t.d.($A$) = sup$\{$t.d.($A/p$) / p prime
ideal of A$\}$).
A tensor
product of AF-domains is perhaps the most natural example of an AF-ring. We
then developed
quite general results for AF-rings, showing that the results do not extend
trivially from integral
domains to rings with zero-divisors.
\bigskip

Our aim in this paper is to extend Wadsworth's results in a different way,
namely to tensor
products of $k-$algebras arising from pullbacks. In order to do this, we
use previous deep
invenstigations of the prime ideal structure of various pullbacks, as in
[1], [2], [3],
[4], [6], [8], [9], [10] and [16]. Moreover, in [14] dimension formulas for
the tensor product of
two particular pullbacks are established and a conjecture on the dimension
formulas for more general
pullbacks is raised; in the present  paper such conjecture is resolved.
\bigskip

Before presenting our main result of section 1, Theorem 1.9, it is
convenient
to recall from [21]
some notation. Let $A \in C$ and let $d, s$ be integers with $0 \leq d \leq
s$. Put
D($s, d, A$) = max$\{$ ht$p[s]$ + min($s, d $ + t.d.($A/p$)) / $p$ prime
ideal of $A\}$. Our
main result is the following : given $R_1 = \varphi_1^{-1}(D_1)$ and $R_2 =
\varphi_2^{-1}(D_2)$ two
pullbacks issued from $T_1$ and $T_2$ respectively. Assume that $D_i$,
$T_i$ are AF-domains
and ht($M_i$) = dim($T_i$), for $i = 1, 2$. Then
\halign{\qquad#&\raggedright#\cr
dim$(R_1 \otimes_k R_2)$ =
        & max$\Big\{$ht$M_1$[t.d.($R_2$)] + D\big(t.d.($D_1$), dim($D_1$),
$R_2\big)$ {\bf ,} \cr
        & \hskip3truecm ht$M_2$[t.d.($R_1$)] + D\big(t.d.($D_2$),
dim($D_2$), $R_1\big)\Big\}$ \cr}

\noindent It turns out ultimately from this
theorem and via a result of Girolami [13] that one may compute (Krull)
dimensions of tensor
products of two $k-$algebras for a large class of (not necessarily
AF-domains) $k-$algebras.
The purpose of Section 2 is to prove the following theorem : with the above
notation,

\noindent\qquad dim$_v(R_1 \otimes_k R_2)$ =
                               min$\Big\{$dim$_vR_1$ + t.d.($R_2$),
dim$_vR_2$ + t.d.($R_1)\Big\}$

\noindent In Section 3 Theorem 3.1 asserts that, with mild restrictions,
tensor products of pullbacks
preserve Jaffard rings. Theorem 3.2 states, under weak assumptions, a
formula similar to that of
Theorem 1.9. It establishes a satisfactory analogue of [4,Theorem 5.4]
(also [1,Proposition 2.7] and
[9, Corollary 1]) for tensor products of pullbacks issued from AF-domains.
We finally focus on
the special case in which $R_1 = R_2$. Some examples illustrate the limits
of our results and
the failure of Wadsworth's results for non-AF-domains.
\vskip1truecm

\noindent{\obf 1. The Krull dimension }
\bigskip

The discussion which follows, concerning basic facts (and notations)
connected with
the prime ideal structure of pullbacks and tensor products of $k-$algebras,
will provide some
background to the main theorem of this section and will be of use in its
proof. Notice first that
we will be concerned with pullbacks (of commutative $k-$algebras) of the
following type :

$$\matrix{R          &\longrightarrow &D\cr
          \downarrow &                &\downarrow\cr
          T          &\longrightarrow &K\cr}$$

\noindent where $T$ is an integral domain with maximal ideal $M$, $K =
T/M$, $\varphi$ is the
canonical surjection from $T$ onto $K$, D is a proper subring of $K$ and $R
= \varphi^{-1}(D)$.
Clearly, $M = (R:T)$ and $D \cong R/M$. Let $p$ be a prime ideal of $R$. If
$M \not\subset p$,
then there is a unique prime ideal $q$ in $T$ such that $q\cap R = p$ and
$T_q = R_p$. However,
if $M \subseteq p$, there is a unique prime ideal $q$ in $D$ such that $p =
\varphi^{-1}(q)$
and the following diagram of canonical homomorphisms

$$\matrix{R_p          &\longrightarrow &D_q\cr
          \downarrow &                &\downarrow\cr
          T_M          &\longrightarrow &K\cr}$$

\noindent is a pullback. Moreover, ht$p$ = ht$M$ + ht$q$ (see [11] for
additional evidence).
We recall from [8] and [1] two well-known results describing how
dimension and
valuative dimension behave under pullback : with the above notation,
dim$R$ = max$\{$dim$T$, dim$D$ + dim$T_M\}$, and
dim$_vR$ = max$\{$dim$_vT$, dim$_vD$ + dim$_vT_M$ + t.d.$(K:D)\}$. However,
while dim$R[n]$
seems not to be effectively computable in general, questions of effective
upper and lower
bounds for dim$R[n]$ were partially answered. The following lower bound
will be useful in the sequel :
dim$R[n]$ $\geq$ dim$D[n]$ + dim$T_M$ + min($n$, t.d.$(K:D)$), where the
equality holds
if $T$ is supposed to be a locally Jaffard domain with ht$M$ = dim$T$ (cf.
[9]). At last, it is a
key  result [13]  that $R$ is an AF-domain  if and only if so are $T$ and
$D$ and t.d.$(K:D) = 0$.
A combination of this result and Theorem 1.9 allows one to compute
dimensions of tensor products of
two $k-$algebras for a large class of (not necessarily AF-domains)
$k-$algebras.
\bigskip

We turn now to tensor products. Let us recall from [21] the following
functions :
let $A, A_1$ and $A_2 \in C$.
Let $p \in$ Spec($A$), $p_1 \in$ Spec($A_1$) and $p_2 \in$ Spec($A_2$).
Let
$d, s$ be integers with
$0 \leq d \leq s$. Set \par

\noindent$\bullet \quad S_{p_1, p_2}$ = $\{ P \in$ Spec$(A_1 \otimes_k
A_2)$ $/ $
$p_1 = P \cap A_1$ and $p_2 = P \cap A_2 \}$ \par

\noindent$\bullet \quad \delta(p_1, p_2)$ = max$\{$ht$P$ $/ $ $P \in
S_{p_1, p_2} \}$ \par

\noindent$\bullet \quad \Delta(s, d, p)$ = ht$p[s]$ + min$(s, d +
t.d.(A/p))$ \par

\noindent$\bullet$\quad D$(s, d, A)$ = max$\{ \Delta(s, d, p)$ $/ $  $p
\in$ Spec($A)\}$ \par

\noindent One can easily check that
dim$(A_1 \otimes_k A_2)$ = max$\{\delta(p_1, p_2)$/ $p_1 \in$ Spec($A_1$)
and $p_2 \in$
Spec($A_2)\}$  (see [21, page 394]).
Let $P \in$ Spec$(A_1 \otimes_k A_2)$ with $p_1 \subseteq P \cap A_1$
and $p_2 \subseteq P \cap A_2$. It is known [21] that $P$ is minimal in
$S_{p_1, p_2}$ if and only
if it is a minimal prime divisor of $p_1 \otimes A_2 + A_1 \otimes p_2$.
This result will be used
to prove a special chain lemma for tensor products of $k-$algebras, which
establishes a somewhat
analogue of the Jaffard's special chain theorem for polynomial rings (see
[7] and [15]).
\bigskip

These facts will be used frequently in the sequel without explicit mention.
\bigskip

The proof of our main theorem requires some preliminaries. The following
two lemmas deal with
properties of polynomial rings over pullbacks, which are probably
well-known, but we have not
located references in the literature.
\bigskip

\noindent {\bf Lemma 1.1.} {\it Let T be an integral domain with maximal
ideal M, K $=$ T/M,
$\varphi$ the canonical surjection from T onto K, D a proper subring of K
and R =
$\varphi^{-1}$(D). Then htp[n] = ht(p[n]/M[n]) + htM[n],
for each positive integer n and each prime ideal p of R such that
M $\subseteq$ p.}
\bigskip

\noindent {\bf Proof.} Since $M \subseteq p$, there is a unique $q \in$
Spec($D$) such that
$p = \varphi^{-1}(q)$ and the following diagram is a pullback

$$\matrix{R_p          &\longrightarrow &D_q\cr
          \downarrow &                &\downarrow\cr
          T_M          &\longrightarrow &K\cr}$$

\noindent By [1, Lemma 2.1 (c)] $MT_M = MR_p$ is a divided prime ideal of
$R_p$.
By [1, Lemma 2.2]
ht$p[n]$ = ht$pR_p[n]$ = ht($pR_p[n]/MR_p[n]$) +  ht$MR_p[n]$ =
ht($p[n]/M[n]$) + ht$M[n]$.
$\diamondsuit$
\bigskip

\noindent {\bf Lemma 1.2.} {\it Let T be an integral domain with maximal
ideal M, K $=$ T/M,
$\varphi$ the canonical surjection from T onto K, D a proper subring of K
and R =
$\varphi^{-1}$(D). Assume T$_M$ and D are locally Jaffard domains.
Then htp[n] = htp + min(n, t.d.(K:D), for each positive integer n and each
prime ideal p of R
such that M $\subseteq$ p.}
\bigskip

\noindent {\bf Proof.} Since $M \subseteq p$, there is a unique $q \in$
Spec($D$) such that
$p = \varphi^{-1}(q)$ and the following diagram is a pullback

$$\matrix{R_p          &\longrightarrow &D_q\cr
          \downarrow &                &\downarrow\cr
          T_M          &\longrightarrow &K\cr}$$

\noindent By [3, Corollary 2.10] ht$p[n]$ = dim$(R_p[n]) - n$. Furthermore,
\halign{\qquad#&\raggedright#\cr
dim($R_p[n]$)  & = ht$M$ + dim($D_q[n]$) + min($n$, t.d.$(K\colon D)$) \cr
               & = ht$M$ + dim$D_q$ + $n$ + min($n$, t.d.$(K\colon D)$) \cr
               & = ht$p$ + $n$ + min($n$, t.d.($K\colon D$)), completing the
proof.
$\diamondsuit$\cr}
\bigskip

The following corollary is an immediate consequence of (1.2) and will be
useful in
the proof of the theorem.
\bigskip

\noindent {\bf Corollary 1.3.} {\it Let T be an integral domain with
maximal ideal M, K $=$
T/M,
$\varphi$ the canonical surjection from T onto K, D a proper subring of K
and R =
$\varphi^{-1}$(D). Assume T$_M$ is a locally Jaffard domain.
Then htM[n] = htM + min(n, t.d.(K:D)), for each positive integer n.}
$\diamondsuit$
\bigskip

We next analyse the heights of ideals of $A_1 \otimes_k A_2$ of the form
$p_1 \otimes_k A_2$,
where $p_1 \in$ Spec($A_1$) and $A_2$ is an integral domain.
\bigskip

\noindent {\bf Lemma 1.4.} {\it Let A$_1$, A$_2$ $\in$ C and $p_1$ be a
prime ideal of $A_1$.
Assume A$_2$ is an integral domain. Then ht($p_1 \otimes_k A_2$) =
ht$p_1[$t.d.$(A_2)]$.}
\bigskip

\noindent {\bf Proof.} Put $t_2$ = t.d.$(A_2)$. Let $Q$ be a minimal prime
divisor of
$p_1 \otimes A_2$ in $A_1 \otimes A_2$. Then $Q$ is minimal in $S_{p_1,
(0)}$, and hence
t.d.$((A_1 \otimes A_2)/Q)$ = t.d.$(A_1 /p_1)$ + $t_2$ by [21, Proposition
2.3]. Furthermore, $Q$
survives in $A_1 \otimes F_2$, where $F_2$ is the quotient field of $A_2$,
whence
ht$Q$ + t.d.$((A_1 \otimes A_2)/Q)$ = $t_2$ + ht$p_1[t_2]$ + t.d.$(A_1
/p_1)$ by [21, Remark 1.b],
completing the proof. $\diamondsuit$
\bigskip

With the further assumption that $A_2$ is an AF-domain, we obtain
the following.
\bigskip

\noindent {\bf Lemma 1.5. (Special chain lemma)} {\it Let A$_1$, A$_2$
$\in$ C and $p_1$ be a prime
ideal of $A_1$.
Assume A$_2$ is an AF-domain. Let $P \in$ Spec$(A_1 \otimes_k A_2)$ such
that $p_1 = P \cap A_1$.
Then htP = ht($p_1 \otimes_k A_2$) + ht(P/($p_1 \otimes_k A_2$)).}
\bigskip

\noindent {\bf Proof.} Since $A_2$ is an AF-domain, by
[21, Remark 1.b] ht$P$ +
t.d.$((A_1 \otimes A_2)/P)$ = $t_2$ + ht$p_1[t_2]$ + t.d.$(A_1 /p_1)$,
where $t_2$ =
t.d.$(A_2)$. A similar argument with $(A_1/p_1) \otimes_k A_2$ in place of
$A_1 \otimes_k A_2$
shows that ht$(P/(p_1 \otimes_k A_2))$ + t.d.$((A_1 \otimes A_2)/P)$ =
$t_2$ + t.d.$(A_1 /p_1)$,
whence ht$P$ = ht$p_1[t_2]$ +  ht$(P/(p_1 \otimes_k A_2))$. The proof is
complete via Lemma 1.4.
$\diamondsuit$
\bigskip

An important case of Lemma 1.5 occurs when $A_2 = k[X_1, . . ., X_n]$ and
hence if $P$ is
a prime ideal of $A_1 \otimes A_2 \cong A_1[X_1, . . ., X_n]$ with $p = P
\cap A_1$, then
ht$P$ = ht$p[n]$ + ht$P/p[n]$. Our special chain lemma may be then viewed
as an analogue
of the Jaffard's special chain theorem (see [7] and [15]). Notice for
convenience that
Jaffard's theorem holds for any (commutative) ring, while here we are
concerned with
$k-$algebras.
\bigskip

To avoid unnecessary repetition, let us fix notation for the rest of this
section and also
for much of section 2 and 3. Data will consist of two pullbacks
of $k-$algebras

$$\matrix{R_1        &\longrightarrow &D_1         &\hskip3truecm
R_2        &\longrightarrow &D_2\cr
          \downarrow &                &\downarrow  &\hskip3truecm
\downarrow &         &\downarrow\cr
          T_1        &\longrightarrow &K_1         &\hskip3truecm
T_2    &\longrightarrow &K_2\cr}$$

\noindent where, for $i = 1, 2$, $T_i$ is an integral domain with maximal
ideal $M_i$,
$K_i = T_i/M_i$, $\varphi_i$ is the canonical surjection from $T_i$ onto
$K_i$, $D_i$ is
a proper subring of $K_i$ and $R_i = \varphi_i^{-1}(D_i)$. Let $d_i$ =
dim$T_i$, $d'_i$ = dim$D_i$,
$t_i$ = t.d.$(T_i)$, $r_i$ = t.d.$(K_i)$ and $s_i$ = t.d.$(D_i)$.
\bigskip

The next result deals with the function $\delta(p_1, p_2)$ according to
inclusion relations
between $p_i$ and $M_i$ ($i = 1, 2$).
\bigskip

\noindent {\bf Lemma 1.6.} {\it Assume $T_1$ and $T_2$ are AF-domains. If
$p_1 \in$ Spec($R_1$) and
$p_2 \in$ Spec($R_2)$ are such that $M_1 \not\subset p_1$ and $M_2
\not\subset p_2$, then

\qquad $\delta(p_1, p_2)$ = min(ht$p_1 + t_2$, $t_1$ + ht$p_2$) $\leq$
min($d_1 + t_2, t_1 + d_2$).}
\bigskip

\noindent {\bf Proof.} By [1, Lemma 2.1 (e)], for $i = 1, 2$, there exists
$q_i \in$ Spec($T_i$)
such that  $p_i = q_i \cap R_i$ and $T_{iq_i}$ = $R_{ip_i}$. So that
$R_{1p_1}$ and $R_{2p_2}$ are
AF-domains, whence $\delta(p_1, p_2)$ = min(ht$p_1 + t_2$, $t_1$ + ht$p_2$)
by [21, Theorem 3.7].
Further, ht$p_1$ $\leq d_1$ and ht$p_2$ $\leq d_2$, completing the proof.
$\diamondsuit$
\bigskip

\noindent {\bf Lemma 1.7.} {\it Assume $T_1$ and $T_2$ are AF-domains. Let
$P \in$ Spec$(R_1 \otimes_k R_2)$, $p_1 = P \cap R_1$ and $p_2 = P \cap
R_2$.
If $M_1 \subseteq p_1$ and $M_2 \not\subset p_2$, then
htP = ht$M_1[t_2]$ + ht$(P/(M_1 \otimes R_2))$}.
\bigskip

\noindent {\bf Proof.} Since $M_2 \not\subset p_2$, $R_{2p_2}$ is
an AF-domain. By Lemma 1.5
ht$P$ = ht$p_1[t_2]$ + ht$(P/(p_1 \otimes R_2))$. Since
$M_1 \subseteq p_1$,
ht$p_1[t_2]$ = ht($p_1[t_2]/M_1[t_2]$) + ht$M_1[t_2]$ by Lemma 1.1.
Hence
\halign{\qquad#&\raggedright#\hfil\cr
ht$P$ &= ht($p_1[t_2]/M_1[t_2]$) + ht$M_1[t_2]$ +
ht$(P/(p_1 \otimes R_2))$\cr
      &= ht$((p_1 \otimes R_2)/(M_1 \otimes R_2))$ + ht$M_1[t_2]$ +
ht$(P/(p_1 \otimes R_2))$\cr
      &$\leq$ ht$M_1[t_2]$ + ht$(P/(M_1 \otimes R_2))$\cr
      &= ht$(M_1 \otimes R_2)$ + ht$(P/(M_1 \otimes R_2))$ \cr
      &$\leq$ ht$P$. $\diamondsuit$\cr}
\bigskip

A similar argument with the roles of $p_1$ and $p_2$ reversed shows that
if $M_1 \not\subset p_1$ and $M_2 \subseteq p_2$, then
htP = ht$M_2[t_1]$ + ht$(P/(R_1 \otimes M_2))$.
\bigskip

Now, we state our last preparatory result, by giving a formula for
dim$((R_1/M_1) \otimes
(R_2/M_2))$ and useful lower bounds for dim$((R_1/M_1) \otimes R_2)$ and
dim$(R_1 \otimes
(R_2/M_2))$.
\bigskip

\noindent {\bf Lemma 1.8.} {\it Assume $T_1$, $T_2$, $D_1$ and $D_2$ are
AF-domains with
dim$T_1$ = ht$M_1$ and dim$T_2$ = ht$M_2$. Then

\noindent {\bf a)} dim$((R_1/M_1) \otimes R_2)$ $\geq$ $d_2$ + min$(s_1,
r_2 - s_2)$ +
                                                             min$(s_1 +
d'_2, d'_1 +  s_2)$.

\noindent {\bf b)} dim$(R_1 \otimes (R_2/M_2))$ $\geq$ $d_1$ + min$(s_2,
r_1 - s_1)$ +
                                                             min$(s_1 +
d'_2, d'_1 +  s_2)$.

\noindent {\bf c)} dim$((R_1/M_1) \otimes (R_2/M_2))$ = min$(s_1 + d'_2,
d'_1 +  s_2)$. }
\bigskip

\noindent {\bf Proof.} a) Since $R_1/M_1 \cong D_1$ is an AF-domain, by
[21, Theorem 3.7]

\qquad dim$((R_1/M_1) \otimes R_2)$ = D$(s_1, d'_1, R_2)$ =
max$\{ \Delta(s_1, d'_1, p_2) / p_2 \in {\rm Spec}(R_2)\}$.

\noindent Let $p_2 \in {\rm Spec}(R_2)$ such that $M_2 \subseteq p_2$.
Then
there is a
unique $q_2 \in$ Spec($D_2$) such that $p_2 = \varphi_2^{-1}(q_2)$ and the
following diagram
is a pullback

$$\matrix{R_{2p_2}          &\longrightarrow  &D_{2q_2}\cr
          \downarrow &                        &\downarrow\cr
          T_{2M_2}          &\longrightarrow  &K_2\cr}$$

\noindent By Lemma 1.2 ht$p_2[s_1]$ = ht$p_2$  + min$(s_1, r_2 - s_2)$.
Since $R_2/p_2$ and
$D_2/q_2$ are isomorphic $k-$algebras, t.d.($R_2/p_2$) = t.d.($D_2/q_2$) =
$s_2$ -  ht$p_2$ +  ht$M_2$, so that

\halign{\qquad#&\raggedright#\hfil\cr
$\Delta(s_1, d'_1, p_2)$ &= ht$p_2[s_1]$ +
min$(s_1, d'_1$, t.d.$(R_2/p_2))$\cr
                         &= ht$p_2$ + min$(s_1, r_2 - s_2)$ +
                                       min$(s_1, d'_1 + s_2 - $
ht$p_2$ + ht$M_2$)\cr
                     &=  min$(s_1, r_2 - s_2)$ + min$(s_1 +  $ ht$p_2 ,
d'_1 + s_2$ + ht$M_2$)\cr
                     &=  ht$M_2$ + min$(s_1, r_2 - s_2)$ + min$(s_1 +  $
ht$q_2 , d'_1 + s_2$)\cr
                     &=  $d_2$ + min$(s_1, r_2 - s_2)$ + min$(s_1 +  $
ht$q_2 , d'_1 + s_2$).\cr}

\noindent b) As in (a) with the roles of $R_1$ and $R_2$ reversed.

\noindent c) It is immediat from [21, Theorem 3.7]. $\diamondsuit$
\bigskip

The facts stated above
provide motivation for setting:\par $\alpha_1 = d_1$ + min$(t_2,  r_1 -
s_1) + d_2$ + min$(s_1,  r_2 - s_2)$ + min$(s_1 + d'_2,
d'_1 + s_2)$;

$\alpha_2 = d_2$ + min$(t_1, r_2 - s_2) + d_1$ + min$(s_2,  r_1 - s_1)$ +
min$(s_1 + d'_2,
d'_1 + s_2)$;

$\alpha_3  = d_1 + d_2$ + min$(r_1, r_2)$ + min$(s_1 + d'_2,  d'_1 +
s_2)$.\par\noindent We shall use
these numbers in the proof of the next theorem and in section 3.
\bigskip

We now are able to state our main result of this section.
\bigskip

\noindent {\bf Theorem 1.9.} {\it Assume $T_1$, $T_2$, $D_1$ and $D_2$
are
AF-domains with
dim$T_1$ = ht$M_1$ and dim$T_2$ = ht$M_2$. Then

\halign{\qquad#&\raggedright#\cr
dim$(R_1 \otimes_k R_2)$ =
     & max$\Big\{$ht$M_1[$t.d.($R_2)]$ + D\big(t.d.($D_1$), dim($D_1$),
$R_2$\big)
                             {\bf ,} \cr
     &\hskip3truecm ht$M_2[$t.d.($R_1)]$ + D\big(t.d.($D_2$),
dim($D_2$), $R_1\big)\Big\}$ \cr}}
\bigskip

\noindent {\bf Proof.} Since
dim$(R_1\otimes R_2 ) \geq$ ht$(M_1\otimes R_2)$ +
dim$((R_1/M_1)\otimes R_2)$,
we have  dim$(R_1\otimes R_2)$ $\geq$ ht$M_1[t_2]$ +
dim$((R_1/M_1)\otimes
R_2)$ by Lemma 1.4.
Similarly,  dim$(R_1\otimes R_2 ) \geq$ ht$M_2[t_1]$ + dim$(R_1\otimes
(R_2/M_1))$. Therefore
it suffices to show that
dim$(R_1\otimes R_2 ) \leq$ max$\{htM_1[t_2]$ + dim$((R_1/M_1)\otimes
R_2)$,   ht$M_2[t_1] $+
dim$(R_1\otimes (R_2/M_2))\}$.

\noindent It is well-known that  dim$(R_1\otimes R_2 )$ =
max$\{\delta
(p_1, p_2)  |
p_1\in$  Spec$(R_1), p_2\in $ Spec$(R_2) \}$.
Let $p_1\in $ Spec$(R_1)$ and $p_2\in$  Spec$(R_2)$.
There are four cases :

\noindent 1.  If  $M_1\not\subset p_1$  and
$M_2\not\subset p_2$, by Lemma
1.6 $\delta (p_1, p_2)$
= min$(htp_1 + t_2,  t_1 + htp_2) \leq \alpha_3 .$

\noindent 2. If $M_1\subseteq  p_1$  and  $M_2\not\subset p_2$,
by Lemma 1.7
$\delta (p_1, p_2) \leq$ ht$M_1[t_2]$ + dim$((R_1/M_1)\otimes R_2)$.

\noindent 3. If  $M_1\not\subset p_1$  and  $M_2\subseteq  p_2$,
by Lemma 1.7
$\delta (p_1, p_2) \leq$ ht$M_2[t_1] $ + dim$(R_1\otimes (R_2/M_2)).$

\noindent 4. If  $M_1\subseteq  p_1$  and  $M_2\subseteq p_2$,
then $\delta
(p_1, p_2) \leq$
max $\{$ ht$M_1[t_2]$ + dim$((R_1/M_1)\otimes R_2)$, ht$M_2[t_1]$ +
dim$(R_1\otimes (R_2/M_2)), \alpha_3\}$. Indeed, put $h =
\delta(p_1,
p_2)$. Pick a chain
$P_0\subset P_1\subset ....\subset P_h$ of $h + 1$ distinct
prime ideals
in $R_1\otimes R_2$ with
$P_h \in S_{p_1, p_2}$. If
$M_1\subset P_0\cap R_1$ and $M_2\subset P_0\cap R_2$, then $h$ =
ht$P_h/P_0 \leq$
dim$((R_1/M_1)\otimes (R_2/M_2)) \leq \alpha_3$. Otherwise,
let  $i$  be
the largest integer
such that  $M_1\not\subset P_i\cap R_1$ and let  $j$ be the
largest integer such that  $M_2\not\subset P_j\cap R_2$. If
$i \not=  j$,
say
$i < j$, by Lemma 1.7
ht$P_j $= ht$M_1[t_2]$+ ht$(P_j/(M_1\otimes R_2))$, whence $h \leq$
ht$M_1[t_2]$ +
ht$(P_h/(M_1\otimes R_2)) \leq$ ht$M_1[t_2]$ +
dim$((R_1/M_1)\otimes R_2)$.
If  $i = j$,
since $M_1\subseteq p_1$, there is a unique  $q_1\in$
Spec$(D_1)$  such that
$p_1 = \varphi^{- 1}_1(q_1)$  and the following diagramm is a
pullback

$$\matrix{R_{1p_1}          &\longrightarrow  &D_{1q_1}\cr
          \downarrow &                        &\downarrow\cr
          T_{1M_1}          &\longrightarrow  &K_1\cr}$$

\noindent Since $M_1\not\subset P_i\cap R_1$,  it follows that
$(P_i\cap
R_1)R_{1p_1}
\subset  M_1 T_{1M_1}   = (R_{1p_1}:T_{1M_1})$ by [1, Lemma 2.1 (c)],
whence
ht$(P_i\cap R_1) \leq$ ht$M_1 - 1 = d_1 - 1$. Similarly,
ht$(P_i\cap R_2) \leq$ ht$M_2 - 1 = d_2 - 1$. Finally,
we get via Lemma 1.6

\halign{\qquad#&\raggedright#\hfil\cr
$h$  &= ht$P_i$ + 1 + ht$(P_h/P_{i + 1})$ \cr
      &$\leq \delta (P_i\cap R_1, P_i\cap R_2) + 1$ +
dim$((R_1/M_1)\otimes
(R_2/M_2)$ \cr
      &= min(ht$(P_i\cap R_1) + t_2,  t_1 +
$ ht$(P_i\cap R_2)) + 1 +$
dim$((R_1/M_1)\otimes
(R_2/M_2)$ \cr
      &$\leq$ min$(d_1 -  1 + t_2,  t_1 + d_2 -  1) +  1 +$
dim$((R_1/M_1)\otimes (R_2/M_2))$ \cr
      &= $\alpha_3$. The fourth case is done.\cr}

Now, let us assume $s_1 \leq r_2 - s_2$. Then

\halign{\qquad#&\raggedright#\hfil\cr
$\alpha_1$  &= $d_1 + $min$(t_2,  r_1 - s_1) + d_2 + s_1 +
$min$(s_1 +
d'_2,  d'_1 + s_2)$ \cr
             &=  $d_1 + $min$(t_2 + s_1,  r_1) + d_2 +
$min$(s_1 + d'_2,
d'_1 + s_2)$ \cr
             &$\geq d_1 + d_2 + $min$(r_1, r_2) +
$min$(s_1 + d'_2,  d'_1 +
s_2)$ \cr
             &= $\alpha_3 $.  \cr}

\noindent If $s_2 \leq r_1 - s_1$, in a similar manner we
obtain $\alpha_2
\geq \alpha_3 $.
Finally,  assume $r_1 - s_1 < s_2$  and  $r_2 - s_2 < s_1$,
so that

\halign{\qquad#&\raggedright#\hfil\cr
$\alpha_1$ &= $\alpha_2$ \cr
            &= $t_1 - s_1 + t_2 - s_2 +
$ min$(s_1 + d'_2,  d'_1 + s_2)$ \cr
         &= min$(t_1 + t_2 - s_2+ d'_2,  t_1 + t_2 - s_1+ d'_1)$ \cr
         &= min(dim$_v  R_1+ t_2,  t_1$ + dim$_v  R_2)$. \cr}

\noindent Hence by [13, Proposition 2.1]

\halign{\qquad#&\raggedright#\hfil\cr
dim$(R_1\otimes R_2 )$  &$\leq$ dim$_v  (R_1\otimes R_2)$ \cr
                        &$\leq$ min(dim$_v R_1 + t_2$, dim$_v R_2
+ t_1)$ \cr
                        &= $\alpha_1 = \alpha_2$ \cr
                        &$\leq$ dim$(R_1\otimes R_2 )$. \cr}

\noindent Finally, one may easily check, via Corollary 1.3 and
Lemma 1.8, that
$\alpha_1 \leq htM_1[t_2] $+ dim$((R_1/M_1)\otimes R_2)$ and
$\alpha_2 \leq$ ht$M_2[t_1] $ + dim$(R_1\otimes (R_2/M_2))$.
$\diamondsuit$
\bigskip

It is still an open problem to compute dim$(R_1\otimes R_2)$ when only
$T_1$ (or $T_2$) is
assumed to be an AF-domain. However, if none of the $T_i$ is an AF-domain
$(i = 1, 2)$, then
the formula of Theorem 1.9 may not hold (see [21, Examples 4.3]).

Now assume $R_i$ is an AF-domain and dim$T_i$ = ht$M_i$ = $d_i$, for each
$i$ = 1, 2.
By [13], $T_i$ and $D_i$ are AF-domains and t.d.($K_i:D_i$) = 0 (that is,
$r_i = s_i$). Further, by [1] dim$R_i$ = dim$T_i$ + dim$D_i$ = $d_i + d'_i$.
Therefore, Theorem 1.9 yields:

\halign{#&\raggedright#\hfil\cr
dim$(R_1\otimes R_2 )$  &= max $\{$ht$M_1[t_2]$ + dim$(D_1\otimes R_2)$,
                                               ht$M_2[t_1]$ +
dim$(R_1\otimes D_2)\}$\cr
                        &= max $\{d_1$ + min(dim$R_2 + s_1, t_2 + d'_1)$,\cr
                        &\hskip5truecm   $d_2$ + min(dim$R_1 + s_2, t_1 +
d'_2)\}$ \cr
                        &= max $\{$min(dim$R_2 + r_1 + d_1, t_2 + d'_1 +
d_1)$,\cr
                        &\hskip5truecm   min(dim$R_1 + r_2 + d_2, t_1 +
d'_2 + d_2)\}$ \cr
                        &= min($t_1$ + dim$R_2$, $t_2$ + dim$R_1$).\cr}

\noindent The upshot is that the formula stated in Theorem 1.9 and
Wadsworth's fomula
match in the particular case where $R_1$ and $R_2$ are AF-domains.

\vskip2truecm

\noindent{\obf 2. The valuative dimension}
\bigskip

It is worth reminding the reader that the valuative dimension
behaves well
with respect
to polynomial rings, that is, dim$_v R[n]$ = dim$_v R$ + n, for each
positive integer n
and for any ring R [15, Theorem 2]. Whereas dim$_v(R_1\otimes R_2)$
seems
not to be effectively
computable in general. In [13] the following useful result is proved:
given
$A_1$ and $A_2$ two
$k-$algebras, then
dim$_v(A_1\otimes A_2)$ $\leq$ min(dim$_vA_1$ + t.d.$(A_2)$,
dim$_vA_2$ + t.d.$(A_1)$).
This section's goal is to compute the valuative dimension for a large
class of tensor products of (not necessarily AF-domains) $k-$algebras. We
are still concerned with those arising
from pullbacks.
\bigskip

The proof of our theorem requires a preliminary result, which
provides a
criterion for
a polynomial ring over a pullback to be an AF-domain.
\bigskip

We first state the following.
\bigskip

\noindent{\bf Lemma 2.1.} {\it Let  A  be an integral domain and n  a
positive integer. Then
$A[n]$ is an AF-domain if and only if, for each prime ideal  p  of A,
ht$p[n]$ + t.d.(A/p)
= t.d.(A).}
\bigskip

\noindent{\bf Proof.} Suppose $A[n]$  is an AF-domain. So for each prime
ideal $p$
of A  ht$p[n]$ + t.d.$(A[n]/p[n])$ = t.d.$(A)$ + $n$ , whence
ht$p[n]$ + t.d.$(A/p)$ = t.d.$(A)$. Conversely, if
$Q\in$ Spec$(A[n])$  and  $p = Q \cap A$, then
by [21, Remark  1.b] ht$Q$
+ t.d.$(A[n]/Q)$ =
$n$ + ht$p[n]$ + t.d.$(A/p)$ since  $A[n] \cong  A\otimes k[n]$.
Therefore,
ht$Q$ + t.d.$(A[n]/Q)$ = $n$ + t.d.$(A)$ = t.d.$(A[n])$.
$\diamondsuit$
\bigskip

\noindent{\bf Proposition 2.2.} {\it Let  T  be an integral domain
with
maximal ideal  M, K = T/M,
and $\varphi$ the canonical surjection. Let  D  be a proper  subring
of  K and
R = $\varphi^{-1}(D )$.
Assume T and D are AF-domains. Let r = t.d.(K)  and  s = t.d.(D).
Then
$R[r - s]$  is an
AF-domain.}
\bigskip

\noindent{\bf Proof.} Let  $p \in$ Spec$(R)$. There are two cases:

\noindent 1.  If $M\not\subset p$, then  $R_p$   is an AF-domain.
So  ht$p$
+ t.d.$(R/p)$ =
t.d.$(R)$.   Further, by [21, Corollary 3.2] ht$p$  = ht$p[r - s]$,
whence
ht$p[r - s]$ +
t.d.$(R/p)$ = t.d.$(R)$.

\noindent 2. If  $M\subseteq  p$ , by Lemma 1.2, ht$p[r - s]$ =
ht$p + r -
s$. Moreover
t.d.$(R/p)$ = $s$ + ht$M -$ ht$p$. Then  ht$p[r - s]$ +
t.d.$(R/p) = r$ + ht$M$ = t.d.$(T)$ = t.d.$(R)$. Consequently, $R[r - s]$ is
an AF-domain
by Lemma 2.1. $\diamondsuit$
\bigskip

     We now present  the main result of this section. We consider two
pullbacks of
k-algebras and use the same notations as in the previous sections.
\bigskip

\noindent{\bf Theorem 2.3.} {\it Let $T_1, T_2, D_1$  and  $D_2$  be
AF-domains, with  dim$T_1$ =
ht$M_1$  and
dim$T_2$ = ht$M_2$, then
dim$_v(R_1\otimes R_2)$ = min(dim$_v R_1 + t_2$ ,
dim$_v R_2 + t_1)$.}
\bigskip

\noindent{\bf Proof.} By Proposition 2.2  $R_1[r_1 - s_1]$  and
$R_2[r_2 -
s_2]$  are AF-domains.
Then $R_1[r_1 - s_1]\otimes R_2[r_2 - s_2]$ is an AF-ring
by [21, Proposition 3.1].
Consequently, by [5, Theorem 2.1]
dim$_v(R_1[r_1 - s_1]\otimes R_2[r_2 - s_2])$ =
dim$(R_1[r_1 - s_1]\otimes
R_2[r_2 - s_2])$

\noindent = min(dim $R_1[r_1 - s_1]$ + t.d.$(R_2[r_2 - s_2])$,
t.d.$(R_1[r_1 - s_1])$ +
         dim $R_2[r_2 - s_2])$

\noindent $\geq$ min$(d_1$ + dim$D_1[r_1 - s_1] + r_1 - s_1 +
t_2 + r_2 - s_2, d_2 + $
       dim$D_2[ r_2 - s_2] + r_1 - s_1 + t_1 + r_2 - s_2)$

\noindent = $r_1 - s_1 + r_2 - s_2 +$ min$( d_1 + d'_1 +
r_1 - s_1 + t_2, d_2 + d'_2 + r_2 - s_2 + t_1)$.

\noindent It turns out that
dim$_v (R_1\otimes R_2) \geq$ min$( d_1 + d'_1 + r_1 -s_1+ t_2,
d_2 + d'_2
+ r_2 - s_2+ t_1)$.
So by [1, Theorem 2.11]
dim$_v(R_1\otimes R_2) \geq$ min(dim$_v R_1$ + $t_2$,  $t_1$  +
dim$_v R_2$).
Therefore by [13, Proposition 2.1] we get
dim$_v(R_1\otimes R_2)$ = min(dim$_v R_1$ + $t_2$ ,
dim$_v R_2$ + $t_1$)
           = $t_1 - s_1 + t_2 - s_2$ + min$(s_1 + d'_2,  d'_1 + s_2)$.
$\diamondsuit$
\vskip1truecm

\noindent{\obf 3. Some applications and examples}
\bigskip

We may now state a stability result. It asserts that, under
mild assumptions on
transcendence degrees, tensor products of pullbacks issued from
AF-domains
preserve
Jaffard rings.
\bigskip

\noindent{\bf Theorem 3.1.} {\it If  $T_1, T_2, D_1$ and $D_2$  are
AF-domains, $M_1$  is
the unique maximal ideal of  $T_1$  with  ht$M_1$ = dim$T_1$  and  $M_2$
is the unique
maximal  ideal of  $T_2$  with
dim$T_2$ = ht$M_2$, then  $R_1\otimes R_2$  is a Jaffard ring if and only
if either
$r_1 - s_1 \leq t_2$  and
$r_2 - s_2 \leq s_1$  or  $r_1 - s_1 \leq s_2$  and  $r_2 - s_2 \leq t_1$}.
\bigskip

\noindent{\bf Proof.} Suppose  $r_1 - s_1 \leq t_2$  and  $r_2 - s_2 \leq
s_1$. Then
$\alpha_1 = t_1 - s_1 + t_2 - s_2 +$ min$(s_1 + d'_2,  d'_1 + s_2) =$
min(dim$_v R_1+ t_2,  t_1
+$  dim$_v R_2$).
By Theorem 1.9 and Theorem 2.3
$\alpha_1 \leq$ dim$ (R_1\otimes R_2) \leq$ dim$_v (R_1\otimes R_2)$ =
min(dim$_v R_1+ t_2,  t_1
+$  dim$_v R_2) = \alpha_1$.
Hence  $R_1\otimes R_2$  is a Jaffard ring. Likewise for  $r_1 - s_1 \leq
s_2$  and  $r_2 - s_2
\leq t_1$. Conversely,
  since  $R_1/M_1 \cong  D_1$  is an AF-domain, by [21, Theorem 3.7]
dim$((R_1/M_1)\otimes R_2)$ = D$(s_1, d'_1, R_2)$ = max$\{\Delta(s_1, d'_1,
p_2)   |   p_2\in$
Spec$(R_2)\}$.
If  $M_2\subseteq  p_2$, by the proof of Lemma 1.8  it follows that
$\Delta(s_1, d'_1, p_2) = d_2 +$ min$(s_1, r_2 - s_2)  +$ min$(s_1 +$
ht$q_2, d'_1 + s_2)$
where  $q_2$  is the unique prime ideal of  $D_2$  such that  $p_2 =
\varphi^{-1}_2(q_2)$.
If $M_2\not\subset  p_2$ , since $R_{2p_2}$ is an AF-domain, then
$\Delta(s_1, d'_1, p_2) =$ ht$p_2[s_1] +$ min$(s_1, d'_1 +$ t.d.$(R_2/p_2))$
             = ht$p_2 +$ min$(s_1, d'_1 +$ t.d.$(R_2/p_2)$
             = min$(s_1 +$ ht$p_2 , d'_1 +$ t.d.$(R_2/p_2)$ + ht$p_2) $
             = min$(s_1 +$ ht$p_2 , d'_1 + t_2)$.
In conclusion, since  ht$p_2 \leq d_2 - 1$ being $M_2$ the unique maximal
ideal of $T_2$ with
dim$T_2$ = ht$M_2$, we get dim$((R_1/M_1)\otimes R_2) =$ max$\{d_2 +$
min$(s_1,  r_2 - s_2) +$
min$(s_1 + d'_2,   d'_1 + s_2)$, min$(s_1 + d_2 - 1 ,  d'_1 + t_2)\}$.
Similarly,
dim$(R_1\otimes (R_2/M_2)) =$ max$\{d_1 +$ min$(s_2,  r_1 - s_1)  +$
min$(s_1 + d'_2,
d'_1 + s_2),$ min$(s_2 + d_1 - 1 ,  d'_2 + t_1)\}$.
Moreover by Theorem 2.3
  dim$_v (R_1\otimes R_2) =$ min(dim$_v R_1 + t_2 ,$  dim$_v R_2 + t_1)$
    = $t_1 - s_1 + t_2 - s_2 +$ min$(s_1 + d'_2,  d'_1 + s_2)$.
Let us assume  $s_1 + d'_2 \leq d'_1 + s_2$. Necessarily, $s_1 + d_2 \leq
t_2 + d'_1$. Applying
Corollary 1.3, we obtain
ht$M_1[t_2]$ + dim$((R_1/M_1)\otimes R_2) = d_1 +$ min$(t_2,  r_1 - s_1) +
d_2+$
min$(s_1,  r_2 - s_2) + s_1 + d'_2$.
On the other hand,  $d_1$ + min$(s_2,  r_1 - s_1) + s_1 + d'_2$ = min$(s_2
+ d_1,  t_1 - s_1)
+   s_1 + d'_2$
= min$(d'_2 +  t_1,  s_2 + d_1 + s_1 + d'_2 ) \geq$ min$(s_2 + d_1 - 1 ,
d'_2 +  t_1)$.
Therefore
ht$M_2[t_1]$ + dim$(R_1\otimes (R_2/M_2)) = d_2 +$ min$(t_1,  r_2 - s_2) +
d_1 +$
min$(s_2,  r_1 - s_1) + s_1 + d'_2$. Consequently,
dim$(R_1\otimes R_2) $= max$\{d_1 +$ min$(t_2,  r_1 - s_1)+ d_2 +$
min$(s_1,  r_2 - s_2) + s_1
+ d'_2, d_2  +$ min$(t_1,  r_2 - s_2) + d_1 +$ min$(s_2,  r_1 - s_1) + s_1
+ d'_2\}$
and
dim$_v (R_1\otimes R_2)$ = $t_1 + t_2 - s_2 + d'_2$ = $d_1 + r_1 + d_2 +
r_2 - s_2 + d'_2$.
Since  $R_1\otimes R_2$  is a Jaffard ring, then either
$d_1 +$ min$(t_2,  r_1 - s_1) + d_2 +$ min$(s_1,  r_2 - s_2) +
s_1 + d'_2$ = $d_1 + r_1
+  d_2 + r_2 - s_2 + d'_2$
or
$d_2 +$ min$(t_1,  r_2 - s_2) + d_1 +$ min$(s_2,  r_1 - s_1) + s_1 + d'_2$
= $d_1 + r_1
+ d_2 + r_2 - s_2 + d'_2 $. Hence  either  $r_1 - s_1 \leq t_2$  and  $r_2
- s_2
\leq s_1$  or  $r_1 - s_1 \leq s_2$  and  $r_2 - s_2 \leq t_1$.
Similar arguments run for  $d'_1 + s_2 \leq s_1 + d'_2$, completing the
proof.$\diamondsuit$
\bigskip

Our next result states, under weak assumptions, a formula similar to that
of Theorem 1.9. It
establishes a satisfactory analogue of [4, Theorem 5.4] (also [1,
Proposition 2.7] and [9, Corollary
1]) for  tensor products of pullbacks issued from AF-domains.
\bigskip

\noindent{\bf Theorem 3.2.}  {\it Assume  $T_1$  and  $T_2$  are
AF-domains, with  dim$T_1$ =
ht$M_1$
and dim$T_2$ = ht$M_2$. Suppose that  either  t.d.$(D_1) \leq$
t.d.$(K_2:D_2)$  or
t.d.$(D_2) \leq$  t.d.$(K_1:D_1)$. Then
dim$(R_1\otimes R_2)$ = max$\{$ht$M_1[t_2]$ + dim$(D_1\otimes R_2)$ ,
ht$M_2[t_1] $
+ dim$(R_1\otimes D_2)\}$.}
\bigskip

Here, since none of $D_i$ is supposed to be an AF-domain $(i = 1, 2)$, the
``dim$(D_i\otimes R_j)$ = D$(s_i, d'_i, R_j)$" assertion is no longer valid
in general
$((i, j) = (1, 2), (2, 1))$. Neither is the ``dim$(D_1 \otimes D_2)$
= min$(s_1 + d'_2, d'_1 +  s_2)$" assertion. Put $\alpha'_3$  =
min$(d_1 + t_2,  t_1 + d_2)$ + dim$(D_1\otimes D_2)$.
\bigskip

\noindent{\bf Proof.} The proof runs parallel with the treatment of Theorem
1.9.
An appropriate modification of its proof yields
dim$(R_1\otimes R_2) \leq$max$\{$ht$M_1[t_2]$ + dim$((R_1/M_1)\otimes R_2)$,
ht$M_2[t_1] +$ dim$(R_1\otimes (R_2/M_2)), \alpha'_3 \}$.
Now there is no loss of generality in assuming that
t.d.$(D_1) \leq$ t.d.$(K_2:D_2)$ (That is, $s_1 \leq r_2 - s_2$).
By Lemma 1.1
and Corollary 1.5
ht$(M_1\otimes R_2)$ + ht$(D_1\otimes M_2)$ = ht$M_1[t_2] $+ ht$M_2[s_1]$ =
ht$M_1$ + min$(t_2,  r_1 - s_1)$ + ht$M_2 $+ min$(s_1,  r_2 - s_2)$ =
min$(d_1 + t_2 + d_2 + s_1,  t_1 + d_2) \geq$ min$(d_1 + t_2 ,  t_1 +
d_2)$. Clearly,
$\alpha'_3$  = min$(d_1 + t_2,  t_1 + d_2)$ + dim$(D_1\otimes D_2) \leq$
ht$(M_1\otimes R_2)$ + ht$(D_1\otimes M_2)$ + dim$(D_1\otimes D_2) \leq$
ht$(M_1\otimes R_2) +$ dim$(D_1\otimes R_2)$. $\diamondsuit$
\bigskip

    We now move to the significant special case in which  $R_1 = R_2$.
\bigskip

\noindent{\bf Corollary 3.3.} {\it Let  T  be an AF-domain with maximal
ideal  M  with
htM = dimT = d, K = T/M, and  $\varphi$ the canonical surjection. Let  D
be a proper subring of  K and  R = $\varphi^{-1}(D)$. Assume $D$ is a
Jaffard domain. Then
dim$(R\otimes R)$ = htM$[t]$ + dim$(D\otimes R)$, where  t = t.d.(T).
If  moreover  t.d.(K:D) $\leq$ t.d.(D), then
dim$(R\otimes R)$ = dim$_v (R\otimes R)$ = t + dim$_v R$.}
\bigskip

\noindent{\bf Proof.} If t.d.$(D) \leq$ t.d.$(K:D)$, the result is
immediate by Theorem 3.2.
Assume t.d.$(K:D) \leq$ t.d.$(D)$. Then
dim$(R\otimes R) \geq$ ht$(M\otimes R)$ + ht$(D\otimes M)$ + dim$(D\otimes
D) \geq$
ht$M[t] +$ ht$M[s] +$ dim$D + t.d.(D)$ = $d  +$ min$(t, $t.d.$(K:D)) + d +$
min$(s, $t.d.$(K:D))$ + dim$D$ + t.d.$(D)$
= min$(t + d,  t - $t.d.$(D)) + d +$ t.d.$(K:D)  +$ dim$D + s = t - s + t +
d' = t +$
dim$_v R \geq$ dim$_v (R\otimes R)$.
This completes the proof.$\diamondsuit$
\bigskip

The following example illustrates the fact that in Theorm 1.9 and
Corollary 3.3
the ``dim$T_i$ = ht$M_i$ $(i = 1, 2)$" hypothesis cannot be deleted.
\bigskip

\noindent{\bf Example 3.4.} Let  $K$  be an algebraic extension
field of  $k$, $T = S^{-1}K[X,Y]$, where
$S = K[X,Y] - ((X)\cup (X - 1,Y))$
and $M = S^{-1}(X)$. Consider the following pullback

$$\matrix{R          &\longrightarrow &k(Y)\cr
          \downarrow &                &\downarrow\cr
          S^{-1}K[X,Y]          &\longrightarrow &K(Y)\cr}$$

\noindent Since  $S^{-1}K[X,Y]$  is an AF-domain and the extension
$k(Y)\subset K(Y)$   is
algebraic, by [13]  $R$  is an AF-domain, so that
dim$(R\otimes R)$ = dim$R$ + t.d.$(R) = 2 + 2 = 4$ by [21,
Corollary 4.2].
However, ht$M[2]$ = ht$M $= 1 and  dim$(k(Y)\otimes R)$ =
min$(2, 1 + 2) = 2$. Hence
ht$M$ + dim$(k(Y)\otimes R)$ = 3. $\diamondsuit$
\bigskip

Theorem 1.9 allows one, via [13], to compute (Krull) dimensions of
tensor
products of two
$k-$algebras for a large class of (not necessarily AF-domains)
$k-$algebras. The next example
illustrate this fact.
\bigskip

\noindent{\bf Example 3.5.} Consider the following pullbacks

$$\matrix{R_1        &\longrightarrow &k(X)         &\hskip3truecm
R_2   &\longrightarrow &k\cr
          \downarrow &           &\downarrow  &\hskip3truecm
\downarrow &     &\downarrow\cr
k(X, Y)[Z]_{(Z)}  &\longrightarrow &k(X, Y) &\hskip3truecm
k(X)[Z]_{(Z)} &\longrightarrow &k(X)\cr}$$

\noindent Clearly, dim$R_1$ = dim$R_2$ = 1 and dim$_v R_1$ =
dim$_v R_2$ =
2. Therefore none
of $R_1$ and $R_2$ is an AF-domain. By Theorem 1.9,  we have
dim$(R_1\otimes R_2)$ = 4. Finally, note that Wadsworth's
formula fails since
min$\{$dim$R_1 $+ t.d.$(R_2)$, dim$R_2$ + t.d.$(R_1)\} = 3$.
$\diamondsuit$
\bigskip

The next example shows that a combination of Theorem 1.9 and Theorem 3.2
allows one to compute
dim$(R_1\otimes R_2)$ for more general $k-$algebras.
\bigskip

\noindent {\bf Example 3.6.} Consider the pullback

$$\matrix{R_1          &\longrightarrow &k\cr
          \downarrow &                &\downarrow\cr
         k(X)[Y]_{(Y)}          &\longrightarrow &k(X)\cr}$$

\noindent $R_1$ is a one-dimensional pseudo-valuation domain with  dim$_v
R_1$ = 2.
Clearly, $R_1$  is not an AF-domain. By Theorem 1.9  dim$(R_1\otimes R_1)$
= 3. Consider now
the pullback
$$\matrix{R_2          &\longrightarrow &R_1\cr
          \downarrow &                  &\downarrow\cr
        k(X, Y, Z)[T]_{(T)}             &\longrightarrow &k(X, Y, Z)\cr}$$

\noindent We have  dim$R_2$ = 2  and  dim$_v R_2$ = 4. The second pullback
does not satisfy
conditions of Theorem 1.9. Applying Theorem 3.2,  we get
dim$(R_1\otimes R_2)$ = max$\{$ht$M_1[4]$ +  dim$(k\otimes~R_2)$,
ht$M_2[2]$ +
dim$(R_1\otimes R_1)\}$ = max$\{2 + 2 ,   2 + 3\}$ = 5. $\diamondsuit$
\bigskip

The next example shows that Corollary 3.3 enables us to construct an
example of an integral
domain $R$  which is not an AF-domain while $R\otimes R$ is a Jaffard ring.
\bigskip

\noindent{\bf Example 3.7.} Consider the pullback

$$\matrix{R        &\longrightarrow &k(X)\cr
          \downarrow &           &\downarrow \cr
k(X, Y)[Z]_{(Z)}  &\longrightarrow &k(X, Y) \cr}$$

\noindent dim$R$ = 1 and dim$_v R$ = 2. Then $R$  is not an AF-domain. By
Corollary 3.3
dim$(R\otimes R)$ = dim$_v(R\otimes R)$ = 5 since t.d.$(k(X,Y): k(X)) <$
t.d.$(R)$.
$\diamondsuit$
\vskip1truecm

\noindent {\obf References}
\bigskip

\item{[1]} D.F. Anderson, A. Bouvier, D.E. Dobbs, M. Fontana and S. Kabbaj,
On Jaffard domains,  Expo. Math. 6 (1988) 145-175.

\item{[2]}  J. T. Arnold, On the dimension theory of overrings of an
integral domain,
Trans. Amer. Math. Soc. 138 (1969) 313-326.

\item{[3]}  J. T. Arnold and R. Gilmer, The dimension sequence of a
commutative   ring,
Amer. J. Math. 96 (1974) 385-408.

\item{[4]}  E. Bastida and R. Gilmer, Overrings and divisorial ideals of
rings of the form D+M,
Michigan Math. J. 20 (1973) 79-95.

\item{[5]} S. Bouchiba, F. Girolami and S. Kabbaj, The dimension of tensor
products of AF-rings,
Lecture  Notes in Pure and Applied  Mathematics, M. Dekker, 185 (1997)
141-154.

\item{[6]} A.Bouvier and S. Kabbaj, Examples of Jaffard domains, Journal of
Pure and Applied Algebra
54 (1988) 155-165.

\item{[7]} J.W. Brewer, P.R. Montgomery, E.A. Rutter and W.J. Heinzer,
Krull dimension of polynomial
rings, Lecture Notes in Mathematics, Springer-Verlag, 311 (1972) 26-45.

\item{[8]} J.W. Brewer and E.A. Rutter, $D$ + $M$ constructions with
general overrings, Michigan
Math. J. 23 (1976) 33-42.

\item{[9]}  P.-J. Cahen, Couples d'anneaux partageant un id\'eal,
Arch. Math. 51 (1988) 505-514.

\item{[10]}  P.-J. Cahen, Construction B, I, D et anneaux localement ou
r\'esiduellement de Jaffard,
Arch. Math. 54 (1990) 125-141.

\item{[11]} M. Fontana, Topologically defined classes of commutative rings,
Ann. Mat. Pura Appl.
123 (1980) 331-355.

\item{[12]}  R. Gilmer, Multiplicative ideal theory, M. Dekker, New York,
1972.

\item{[13]}  F. Girolami, AF-rings and locally Jaffard rings,
Lecture  Notes in Pure and Applied  Mathematics, M. Dekker, 153 (1994)
151-161.

\item{[14]}  F. Girolami and S. Kabbaj, The dimension of the tensor
product
of two particular
pullbacks, Proc. Padova Conference "Abelian groups and modules", Kluwer
Academic Publishers,
(1995) 221-226.

\item{[15]}  P. Jaffard, Th\'eorie de la dimension dans les anneaux de
polyn\^omes,
M\'em. Sc. Math., 146, Gauthier-Villars, Paris, 1960.

\item{[16]} S. Kabbaj, Quelques probl\`emes sur la th\'eorie des spectres
en alg\`ebre commutative,
Thesis of Doctorat Es-Sciences, University of F\`es, F\`es, Morocco, 1989.

\item{[17]} I. Kaplansky, Commutative Rings, University of Chicago Press,
Chicago, 1974.

\item{[18]}  H. Matsumura, Commutative ring theory, Cambridge University
Press, Cambridge, 1989.

\item{[19]}  M. Nagata, Local rings, Interscience,  New York, 1962.

\item{[20]}  R.Y. Sharp, The dimension of the tensor product of two field
extensions,  Bull. London Math. Soc. 9 (1977) 42-48.

\item{[21]} A.R. Wadsworth, The Krull dimension of tensor products of
commutative algebras over
a field,  J. London Math. Soc. 19 (1979) 391-401.

\end